\documentclass[12pt]{article}

\usepackage{amsfonts,amsmath,amssymb}

\numberwithin{equation}{section}

\title{New Solutions of $d=2x^3+y^3+z^3$ }

\author{Allan J. MacLeod\\Statistics, O.R. and Mathematics Group,\\School of Science,\\
University of the West of Scotland,\\High St.,  Paisley,\\Scotland.  PA1 2BE\\
(e-mail: allan.macleod@uws.ac.uk) }

\date{}

\begin{document}

\maketitle

\begin{abstract}
{\noindent
We discuss finding large integer solutions of $d=2x^3+y^3+z^3$ by using Elsenhans and Jahnel's
adaptation of Elkies' LLL-reduction
method. We find $28$ first solutions for $|d|<10000$.
}

\vspace{0.5cm}

\end{abstract}

\section{Introduction}
Fermat's Last Theorem states that
\begin{displaymath}
x^3+y^3=z^3
\end{displaymath}
has no solutions in non-zero integers.

Over the years, there has been interest in small values of
\begin{equation}
d=x^3+y^3-z^3=x^3+y^3+(-z)^3
\end{equation}

Since $u^3 \equiv 0,1,8 (\bmod \, 9)$ for all integers $u$, it is an easy enumeration of possibilities
to see that {\bf NO} solutions can exist for $d \equiv 4,5 (\bmod \, 9)$. It is conjectured that solutions
exist for all other values of $d$ and it is an on-going project to find such solutions for small $d$.

For
\begin{equation}
d=2x^3+y^3+z^3
\end{equation}
there are no congruence restrictions on $d$. The problem is discussed in Guy \cite{guy}, and Koyama
presents a specialized search method in \cite{koy} but
much less work seems to have been done on this problem. It appears
somewhat easier than $(1.1)$ as there are only $3$ values of $|d|<1000$ for which no solution is known,
namely $d=148,671,788$, and only a total of $62$ values for $|d|<10000$, according to Hisanori Mishima's web-site
\cite{mish}.

Recent, Elsenhans and Jahnel \cite{elsj} used Elkies' method \cite{elk}
to discover previously unknown first solutions to $(1.1)$ for some $|d|<1000$, but with $x$, $y$ and
$z$ very large in size. The paper \cite{elsj} does not contain many technical details. Jahnel's web-site
gives a zipped file containing the C-code used, together with a small Word document describing the
method. Unfortunately, all the comments in the code and most of the Word file are in German, which the
author last seriously studied in the 1970s.

In an attempt to reverse-engineer their methodology, equation $(1.2)$ was used as a
test-case. Since this problem might be popular with amateur mathematicians, the following
description is intended to be as simple as possible.

\section{Basic Theory}
We are considering solutions which cannot be found by simple searching over small values of $x,y,z$.
We assume $|d|<10000$, and $\min(|y|,|z|)>100$, so that there must be a mixture of positive and negative
values of $x,y,z$.

We write the problem as
\begin{equation}
2x^3+y^3=z^3+d
\end{equation}
with $d$ positive or negative, and, thus, we can assume, without loss of generality, that $z>0$.
Then at least one of $x,y$ must be positive.

We have,
\begin{equation}
\left( \frac{y}{z} \right)^3+2\left(\frac{x}{z}\right)^3=1+\frac{d}{z^3}
\end{equation}
and, if $z$ is large, we are effectively looking for rational points close to the curve
\begin{equation}
Y^3=1-2X^3
\end{equation}
and it is for this type of problem that Elkies initially proposed his method.

We have $3$ intervals of interest
\begin{enumerate}
\item[(a)] If $X>0,Y>0$, then $0<X<\sqrt[3]{0.5}$.
\item[(b)] If $Y<0$, then $X>\sqrt[3]{0.5}$.
\item[(c)] For any $X<0$ we have $Y>0$.
\end{enumerate}

Define
\begin{equation}
Y=f(X)=(1-2X^3)^{1/3}
\end{equation}
and consider, for ease of exposition, $0\le X\le \sqrt[3]{0.5} \approx 0.7937$.

We have,
\begin{displaymath}
\frac{dY}{dX}=\frac{-2X^2}{(1-2X^3)^{2/3}} \hspace{2cm} \frac{d^2Y}{dX^2}=\frac{-4X}{(1-2X^3)^{5/3}}
\end{displaymath}

Let $X_0$ be a random point in the interval $(0,0.7937)$, with $Y_0=f(X_0)$ so that
\begin{displaymath}
\frac{dY}{dX}=\frac{-2X_0^2}{Y_0^2} \hspace{3cm} \frac{d^2Y}{dX^2}=\frac{-4X_0}{Y_0^5}
\end{displaymath}

Consider the interval $I=[X_0-h/2,X_0+h/2]$ with $h$ strictly positive, where we usually assume the
end-points stay inside the basic interval for $X$. In this interval, we wish to surround the curve
by a parallelogram. Elsenhans and Jahnel use the German word {\it fliesen} $\equiv$ flagstone for
such a parallelogram. We assume the parallelogram has $2$ edges parallel to the $Y$-axis, but the
other $2$ edges are angled downwards. In fact, we only need to specify the equation of the middle
downward line of the parallelogram.

We know the gradient of the curve at $(X_0,Y_0)$ is $A=-2X_0^2/Y_0^2$, so we look for sides of the
parallelogram with gradient $A$. The equation of the tangent itself is
\begin{equation}
Y=A\,X+(Y_0-A\,X_0)
\end{equation}
but the properties of the curve show that this line is always above the curve.

What we want is a line
\begin{displaymath}
Y=A\,X+B
\end{displaymath}
with $B<Y_0-A\,X_0$ which lies both above and below the curve in $I$.

We have, using Taylor series,
\begin{equation}
f(X_0+\frac{h}{2})\approxeq f(X_0)+f'(X_0)\frac{h}{2}+f''(X_0)\frac{h^2}{8}
\end{equation}
\begin{displaymath}
=Y_0+A\frac{h}{2}+f''(X_0)\frac{h^2}{8}
\end{displaymath}
which is just
\begin{displaymath}
A\,(X_0+h/2)+(Y_0-A\,X_0+h^2f''(X_0)/8)
\end{displaymath}

If we set $B=Y_0-A\,X_0+h^2f''(X_0)/8$, the resulting line lies totally below the curve,
remembering that the second derivative is negative, but if we define
\begin{displaymath}
B=Y_0-A\,X_0+h^2f''(X_0)/16
\end{displaymath}
we get a line which satisfies the requirement of above and below the curve.

Using the above discussion, we look for rational $X=x/z,Y=y/z$ with
\begin{equation}
| X-X_0|<h/2
\end{equation}
and
\begin{equation}
|Y-(AX+B)|<K
\end{equation}
where $K>0$ but usually very small.

If we restrict $z$ to $0<z<L$, this gives
\begin{equation}
|x-X_0z|<h|z|/2<hL/2
\end{equation}
and
\begin{equation}
|y-A\,x-B\,z|<K|z|<KL
\end{equation}

We, thus, have the three restrictions
\begin{displaymath}
\left|\frac{2}{hL}x-\frac{2X_0}{hL}z\right|<1
\end{displaymath}
and
\begin{displaymath}
\left|\frac{1}{KL}y-\frac{A}{KL}x-\frac{B}{KL}z\right|<1
\end{displaymath}
and
\begin{displaymath}
0 < \frac{z}{L} < 1
\end{displaymath}

Define the matrix
\begin{displaymath}
F = \left( \begin{array}{ccccc}
\frac{2}{hL}&\,&0&\,&\frac{-2X_0}{hL}\\
\,&\,&\,&\,&\,\\
\frac{-A}{KL}&\,&\frac{1}{KL}&\,&\frac{-B}{KL}\\
\,&\,&\,&\,&\,\\
0&\,&0&\,&\frac{1}{L}
\end{array} \right)
\end{displaymath}
with the columns denoted $\underline{f}_1$, $\underline{f}_2$ and $\underline{f}_3$.

Then we are looking for vectors in the lattice
\begin{displaymath}
\underline{v}=x\underline{f}_1+y\underline{f}_2+z\underline{f}_3 \hspace{1cm} x,y,z \in \mathbb{Z}
\end{displaymath}
which have small components.

This is what the LLL algorithm does, see Cohen \cite{coh}. So we apply LLL-reduction to $F$,
and, hopefully, find a solution to the problem with $|d|$ small.

As an example, let $h=0.001$, $K=0.00001$ and $L=1000$. $X_0=0.31415$ gives $A=-0.205986$ and
$B=1.043599$, so that we have
\begin{displaymath}
F = \left( \begin{array}{rrr}
2&0&-0.6283\\20.5986&100&-104.3599\\0&0&0.001
\end{array} \right)
\end{displaymath}

Applying the function {\bf qflll} from the software package Pari version 2.3.4 gives the integer
matrix
\begin{displaymath}
M=\left( \begin{array}{rrr}
-15&-74&-313\\-47&-230&-976\\-48&-235&-997
\end{array} \right)
\end{displaymath}
and we have
\begin{displaymath}
F\,M=\left( \begin{array}{rrr}
0.1584&-0.3495&0.4151\\0.2962&0.2801&-0.5415\\-0.0480&-0.2350&-0.9970
\end{array} \right)
\end{displaymath}

Investigating the first column of $M$ quickly gives
\begin{displaymath}
2(-15)^3+(-47)^3+48^3=19
\end{displaymath}
whilst the second column gives
\begin{displaymath}
2(-74)^3+(-230)^3+235^3=427
\end{displaymath}

For the other two intervals $X>\sqrt[3]{0.5}$ and $X<0$, the only difference is the behavior
of the derivatives, but the resulting formulae all turn out the same.

\section{Pyramid Searching}
The first version of the algorithm, which we implemented in Pari, generated $F$,
computed $M$ and just investigated the $3$ columns separately as possible solutions. The columns, however,
give small vectors in the lattice but not necessarily linked to small values of $|d|$.

To ensure we do not miss such solutions we need to investigate further. As pointed out by
Elsenhans and Jahnel, the relations $(2.9)$ and $(2.10)$ give a feasible region which is, in fact,
a pyramid. The vertex is at the origin, and the four other vertices at the solutions of
\begin{equation}
F \, \left( \begin{array}{r} x\\y\\z \end{array} \right) =
\left( \begin{array}{r} \pm 1\\ \pm1\\1 \end{array} \right)
\end{equation}

A look at the numerical solution in the previous sections shows that there can be a huge
number of lattice points within this pyramid, so that testing them all will be impractical.
The matrix $M$, however, is unimodular (determinant $=\pm 1$), so
\begin{equation}
F \, \left( \begin{array}{r} x\\y\\z \end{array} \right) =
F\,M \, \left( \begin{array}{r} x'\\y'\\z' \end{array} \right) =
\left( \begin{array}{r} \pm 1\\ \pm1\\1 \end{array} \right)
\end{equation}
with $x',y',z' \in \mathbb{Z}$.
The matrix $H=F\,M$ is much nicer to deal with and gives a much smaller pyramid to search.

We work out the equations of the $8$ edges of the pyramid, and loop through the possible $z'$
values. For each $z'=$constant plane, we work out the $(x',y')$ limits and loop through these.
This does not work for $z'=0$ so we specify the $(x',y',0)$ triples. We, especially, do not want
to miss the first column of $M$ as it is, very often, the smallest non-zero vector in the lattice.

\section{Practical Implementation}
If we study the numerical results in Koyama, we find that the solutions all have $X$ roughly
in the range $[-10,10]$, so we generate random values of $X_0$ in such a range. For input
values of $X_0,h,K,L$ a Pari $2.3.4.$ program generates $F$ and the Pari code
{\bf qflll} generates the LLL-reduction matrix $M$ and, hence, $F\,M$. We then use the pyramid method
of section $3$ to search for a solution
amongst the $62$ unsolved values of $d$ from Mishima's web pages.

A large amount of experimental testing suggests that $K=O(h^2)$ and $L=O(1/h)$ are good choices,
but these are fairly imprecise relationships.

The program is very small and was run on several lap-tops. The results found so far are

\begin{center}
TABLE 1\\Solutions of $d=2x^3+y^3+z^3$\\
\begin{tabular}{rrrrrrr}
$\;$&$\;$&$\;$&$\;$&$\;$&$\,$&$\,$ \\
$d\,\,$&$\,$&$x\,\,\,\,\,\,$&$\,$&$y\,\,\,\,\,\,$&$\,$&$z\,\,\,\,\,\,$\\
1247&$\,$&26478194&$\,$&108525095&$\,$&-109565866\\
1462&$\,$&5032942&$\,$&24350809&$\,$&-24493307\\
1588&$\,$&6273700&$\,$&-6232583&$\,$&-6314285\\
2246&$\,$&5775101&$\,$&-2344975&$\,$&-7194061\\
2822&$\,$&-6205213&$\,$&-3630235&$\,$&8070731\\
3307&$\,$&13896334&$\,$&-5049077&$\,$&-17367182\\
3335&$\,$&-96695533&$\,$&-16135834&$\,$&121923017\\
3641&$\,$&-30068863&$\,$&-43454404&$\,$&51479399\\
4990&$\,$&-282120164&$\,$&-222569165&$\,$&382436587\\
5188&$\,$&274919617&$\,$&-240828017&$\,$&-302168405\\
5279&$\,$&-77345653&$\,$&-372742111&$\,$&374949254\\
5620&$\,$&4200208&$\,$&3637291&$\,$&-5811935\\
5629&$\,$&-106122938&$\,$&-172761995&$\,$&196148872\\
6707&$\,$&-8265969&$\,$&4019725&$\,$&10210900\\
6980&$\,$&-22045411&$\,$&-12306889&$\,$&28558571\\
7097&$\,$&-359099686&$\,$&-179826055&$\,$&461715044\\
7177&$\,$&-5830313&$\,$&-12164090&$\,$&12998491\\
7323&$\,$&-5182681&$\,$&-5155563&$\,$&7461728\\
7519&$\,$&-13121624&$\,$&8938207&$\,$&15610924\\
7853&$\,$&-20412862&$\,$&10234469&$\,$&25166600\\
8114&$\,$&-5609033023&$\,$&-1349280025&$\,$&7083296297\\
8380&$\,$&-9992594&$\,$&-52534419&$\,$&52774343\\
8572&$\,$&-9632054&$\,$&-16002605&$\,$&18054625\\
8644&$\,$&127473520&$\,$&524520025&$\,$&-529492061\\
8887&$\,$&-8808071&$\,$&-17339003&$\,$&18738346\\
9274&$\,$&-17639675&$\,$&2180206&$\,$&22217602\\
9589&$\,$&-27316976&$\,$&-79014746&$\,$&81134053\\
9850&$\,$&-874953287&$\,$&582754948&$\,$&1045170154
\end{tabular}
\end{center}

\end{document}